\documentstyle[12pt]{article}
\newcommand{\const}{\mathop{\rm const}\limits}

\newcommand{\supp}{\mathop{\rm supp}\limits}

\newcommand{\mes}{\mathop{\rm mes}\limits}

\begin{document}

\begin{center}

\vspace{3mm}

{\bf Sharp Estimates for Module of Continuity  }\par

\vspace{4mm}

{\bf of Fractional Integrals and Derivatives. } \\

\vspace{4mm}

 $ {\bf E.Ostrovsky^a, \ \ L.Sirota^b } $ \\

\vspace{4mm}

$ ^a $ Corresponding Author. Department of Mathematics and computer science, Bar-Ilan University, 84105, Ramat Gan, Israel.\\

E-mail: eugostrovsky@list.ru\\

\vspace{3mm}

$ ^b $  Department of Mathematics and computer science. Bar-Ilan University,
84105, Ramat Gan, Israel.\\

E-mail: sirota3@bezeqint.net \\

\vspace{4mm}
                    {\sc Abstract.}\\

 \end{center}

 \vspace{3mm}

  We derive the bilateral estimates for the module of continuity of
 the fractional integrals and derivatives for the functions  from the classical Lebesgue-Riesz spaces. \par

  \vspace{4mm}

 \par

\vspace{3mm}


\vspace{4mm}

\section{Notations. Statement of problem.}

\vspace{4mm}

 "Fractional derivatives have been around for centuries  but recently they have
found new applications in physics, hydrology and finance", see  \cite{Meerschaert1}. \par
 Another applications: in the theory of Differential Equations are described in \cite{Miller1};
in statistics see in \cite{Bapna1}, \cite{Ostrovsky8}; see also \cite{Golubev1}, \cite{Enikeeva1}; in the theory of
integral equations etc. see in  the classical monograph \cite{Samko1}.  \par

\vspace{4mm}

  Let $ \alpha=\const \in (0,1); $ and let $ g=g(x), \  x \in R_+ $  be measurable numerical function. The fractional derivative
 of a Riemann-Liouville type of order $  \alpha: \  D^{\alpha}[g](x)  = g^{(\alpha)}(x) $  \cite{Riemann1}, \cite{Liouville1}
 is defined as follows: $ \Gamma(1-\alpha)  g^{(\alpha)}(x) =  $

$$
\Gamma(1-\alpha) \ D^{\alpha}[g](x)  =
\Gamma(1-\alpha) \ D^{\alpha}_x[g](x) \stackrel{def}{= } \frac{d}{dx} \int_0^x \frac{g(t) \ dt}{(x-t)^{\alpha}}. \eqno(1.1)
$$
see, e.g. the classical monograph of S.G.Samko, A.A.Kilbas and O.I.Marichev \cite{Samko1}, pp. 33-38; see also \cite{Miller1}.\par

 The case when $ \alpha \in (k, k + 1), \ k = 1,2,\ldots $ may be considered analogously, through the suitable derivatives of integer order. \par

 Hereafter $ \Gamma(\cdot) $ denotes the ordinary Gamma  function. \par

 We agree to take $ D^{\alpha}[g](x_0) = 0,  $ if at the point $  x_0 $ the expression $  D^{\alpha}[g](x_0)  $ does not exists. \par

Notice that the operator of the fractional derivative is non-local, if $ \alpha $ is not integer non-negative number.\par

 Recall also that the fractional  integral $  I^{(\alpha)}[g](x) = I^{\alpha}[g](x)  $ of a Riemann-Liouville type of
an order $ \alpha, 0 < \alpha < 1 $ is defined as follows:

$$
I^{(\alpha)}[g](x) \stackrel{def}{=} \frac{1}{\Gamma(\alpha)} \cdot \int_0^x \frac{g(t) \ dt}{(x-t)^{1 - \alpha}}, \ x,t > 0. \eqno(1.2)
$$
 It is known  (theorem of Abel, see \cite{Samko1}, chapter 2, section 2.1)
 that the operator $ I^{(\alpha)}[\cdot]  $ is inverse to the fractional derivative operator $ D^{(\alpha)}[\cdot],  $
at least  in the class of absolutely continuous functions. \par

 Another approach to the introducing of the fractional derivative, more exactly, the fractional Laplace operator
leads us to the using of Fourier transform

$$
 F[f](t) = \int_{R^d} e^{i (t,x) } \ f(x) \ dx
$$
in the space $  R^d, \ d = 1,2,\ldots:  $

$$
R_{\alpha,F}[f] := C_1(\alpha,d) F^{-1} \left[ |x|^{\alpha} F[f](x)  \right], \ 0 < \alpha < d,
$$
which leads us in turn up to multiplicative constant to the well-known Riesz potential

$$
R_{\alpha}[f] \stackrel{def}{=} \int_{R^d} \frac{f(y) \ dy}{ |x-y|^{d - \alpha}}, \ 0 < \alpha < d. \eqno(1.3)
$$

 Hereafter $ (t,x) = \sum_{m=1}^d t_m x_m, \ |x| = \sqrt{ (x,x) }, \ t,x \in R^d.  $ \par

\vspace{4mm}

 {\bf  We intent to obtain in this short article  the bilateral estimates  for module of continuity of
 the fractional integrals and derivatives for the functions  from the classical Lebesgue-Riesz spaces.} \par

\vspace{3mm}

 The  classical Lebesgue-Riesz  $  L_p $ estimations for the fractional integrals and derivatives are investigated in many works, see, e.g.
  \cite{Adams1}, chapter 2,3;  \cite{Adams2}, \cite{Fuglede1}, \cite{Gatto1},  \cite{Hardy0}, \cite{Hardy1}, \cite{Hardy2}, \cite{Hardy3},
 \cite{Lieb3}, \cite{Harboure1}, \cite{Lieb3}, \cite{Leoni1}, \cite{Maly1}, \cite{Ostrovsky10},\cite{Ostrovsky11}, \cite{Ostrovsky12},
 \cite{Ostrovsky13}, \cite{Ostrovsky19}, \cite{Samko1}, chapters 2-3 etc.  The module of continuous estimates ones (previous works) see
 in \cite{Hajibayov1}, \cite{Mizuta1} - \cite{Mizuta4}, \cite{Samko1}, pp, 66-71. \par

  Note that in the articles \cite{Hajibayov1}, \cite{Mizuta2}, \cite{Mizuta4} is considered the case of the so-called Lebesgue-Riesz
 spaces with variable exponent $ p = p(x) $ for the function $  f(\cdot). $  \par

  Recall that the classical Lebesgue-Riesz  $  L(p) $ norm $  |f|_p  $ of a function $  f $ is defined by a formula

$$
|f|_p = \left[ \int_{R^d} |f(x)|^p \ dx  \right]^{1/p}, \ 1 \le p < \infty
$$
or correspondingly

$$
|f|_p = \left[ \int_{R_+} |f(x)|^p \ dx  \right]^{1/p}, \ 1 \le p < \infty.
$$

 \vspace{4mm}

 \section{Module of continuity of the fractional derivatives.}

 \vspace{3mm}

 The module of continuity $ \omega(f,h), \ h \ge 0  $ of (uniformly continuous) function $  f: R^d \to R,  $ or $ f:(0,b) \to R,
 \ b = \const \in (0,\infty]  $ is defined as usually as follows

$$
\omega(f,h) = \sup_{(x,y): |x - y| \le h } |f(x) - f(y)|. \eqno(2.1)
$$
 Let $  f: (0,b) \to R $ be uniformly continuous function such that $  f(0) = 0 $  and let $ \alpha \in (0,1). $
 The inequality

$$
\omega(D^{\alpha}[f], h)   \le C_D(\alpha) \ \int_0^h \frac{\omega(f,t) \ dt}{t^{1 + \alpha}} \eqno(2.2)
$$
is proved, e.g. in \cite{Samko1}, p. 250-253, theorem 3.16. \par

 Define for simplicity the following  set of continuous functions

$$
S(\alpha) \stackrel{def}{=} \left\{ f: f(0) = 0,  \int_0^h \frac{\omega(f,t) \ dt}{t^{1 + \alpha}} < \infty \right\}. \eqno(2.3)
$$

We define the minimal value os the "constant" $  C_D(\alpha) $ as $  K_D(\alpha): $

$$
 K_D(\alpha) \stackrel{def}{=} \sup_{h \in (0,1)} \ \sup_{f \ne \const, f \in S(\alpha) }
\left[  \omega(D^{\alpha}[f], h): \ \int_0^h \frac{\omega(f,t) \ dt}{t^{1 + \alpha}} \right]. \eqno(2.4)
$$

 It is proved in fact in  \cite{Samko1}, p. 250-253

$$
K_D(\alpha) \le  C \alpha^{-1} \Gamma(1 - \alpha), \ \alpha \in (0,1). \eqno(2.5)
$$
where $  C $ is an absolute constant.\par

\vspace{3mm}

{\bf Proposition 2.1.}

$$
K_D(\alpha) \ge   \Gamma(1 - \alpha), \ \alpha \in (0,1). \eqno(2.6)
$$

\vspace{3mm}

{\bf Proof.} We set $ b = 1 $ and choose $ g(x) = x^{\beta}, \ \beta = \const \in (\alpha, 1], $ (test function), then
$  g(\cdot) \in S(\alpha) $ and $ \omega(g,h) = g(h) = h^{\beta}. $ Further,

$$
\int_0^h \frac{\omega(g,t)}{t^{1 + \alpha}} \ dt = \frac{h^{\beta - \alpha}}{\beta - \alpha},
$$

$$
D^{\alpha}[g] = \frac{d}{dx} x^{1 + \beta - \alpha} \ \int_0^1 z^{\beta} \ (1 - z)^{-\alpha} \ dz  =
$$

$$
(1 + \beta - \alpha) \ x^{\beta - \alpha} \  B(\beta + 1, 1 - \alpha), \eqno(2.7)
$$
where as ordinary $ B(\cdot, \cdot) $ denotes the usually Beta function;

$$
\omega(D^{\alpha}[g],h) = (1 + \beta - \alpha) \ h^{\beta - \alpha} \  B(\beta + 1, 1 - \alpha),
$$

$$
\omega(D^{\alpha}[g],h): \int_0^h \frac{\omega(g,t)}{t^{1 + \alpha}} \ dt =
(1 + \beta - \alpha) (\beta - \alpha) B(\beta + 1, 1 - \alpha),
$$
and we deduce using the well-known recursion for the Gamma function

$$
K_D(\alpha) \ge \frac{\Gamma^2(1 - \alpha) \ \Gamma(1 + \beta)}{\Gamma(\beta - \alpha)}, \eqno(2.8)
$$
therefore

$$
K_D(\alpha) \ge \sup_{\beta \in (\alpha,1)} \left[ \frac{\Gamma^2(1 - \alpha) \ \Gamma(1 + \beta)}{\Gamma(\beta - \alpha)} \right] =
$$

$$
\Gamma^2(1 - \alpha) \sup_{\beta \in (\alpha,1)} \left[ \frac{\Gamma(\beta + 1)}{\Gamma(\beta - \alpha) } \right] = \Gamma(1 - \alpha).
\eqno(2.9)
$$

\vspace{4mm}

 \section{Module of continuity of the fractional integrals.}

\vspace{3mm}

 Let again $  \alpha \in (0,1) $ and let now $  p > 1/\alpha, \ f \in L_p(0, \infty); $
then the function  $ g(x) := I^{\alpha}[f](x) $ is continuous and moreover

$$
 \omega(g, h) = o(h^{\alpha - 1/p}), \ h \to 0+,
$$
see \cite{Samko1}, pp. 104-110. \par

 Denote

$$
Z(\alpha,p) = \left[ \frac{p-1}{\alpha p - 1} \right]^{1 - 1/p}, \eqno(3.1)
$$
and for any function $ f \in L_p(0, \infty), \  $

$$
\Delta_p(f,h) \stackrel{def}{=} \sup_{|\delta| < h} \sup_x \left[ \int_x^{x + \delta} |f(t)|^p \ dt  \right]^{1/p}, \eqno(3.2)
$$
where the function $  f(t) $ is presumed to be continued as zero for negative values $  t. $
Evidently,
$$
\lim_{h \to 0+} \Delta_p(f,h) = 0, \ f \in L_p,
$$
and  \hspace{4mm} $ \Delta_p(f,h) \le |f|_p. $ \par
 We conclude after some calculations  following \cite{Samko1}, pp, 66-71

 $$
 \Gamma(\alpha) |I^{\alpha}[f](x)| \le Z(\alpha,p) \cdot x^{\alpha - 1/p} \cdot \Delta_p(f,x), \ f \in L_p, \ x \in [0,1]; \eqno(3.3)
 $$

$$
\omega ( \Gamma(\alpha) I^{\alpha}[f], h) \le 4 \cdot Z(\alpha,p) \cdot h^{\alpha - 1/p} \cdot \Delta_p(f,h), \ f \in L_p, \ h \in [0,1]; \eqno(3.4)
$$
see also  \cite{Ostrovsky61} for technical details. \par
 In turn,

$$
 \Gamma(\alpha) |I^{\alpha}[f](x)| \le Z(\alpha,p) \cdot x^{\alpha - 1/p} \cdot |f|_p, \ f \in L_p, \ x \in [0,1]; \eqno(3.3a)
 $$

$$
\omega ( \Gamma(\alpha) I^{\alpha}[f], h) \le 4 \cdot Z(\alpha,p) \cdot h^{\alpha - 1/p} \cdot |f|_p, \ f \in L_p, \ h \in [0,1].\eqno(3.4a)
$$

 Recall that here $ \alpha > 1/p; $ the case $  \alpha = 1/p $ is considered  in \cite{Samko1}, p. 69. The case of greatest values $  \alpha $
can by reduced  to the considered here by means of differentiating,  see \cite{Samko1}, p. 69-77.\par

\vspace{3mm}

 We intent to show in this section that the exponent $  \alpha - 1/p $ is exact.\par

\vspace{3mm}

{\bf A. \ Case } $ X = \{x\} = (0,1). $ \par

\vspace{3mm}

{\bf Proposition  3.A.} For all the values $  \epsilon \in (0, (\alpha - 1/p)/2) $  there exists a function $  f_0, \ f_0 \in L_p(0,1),$
for which

$$
\omega(I^{\alpha}[f_0], h) \ge C(\alpha,p) \cdot h^{\alpha - 1/p + \epsilon}, \ C(\alpha,p) > 0, \ h \in (0,1). \eqno(3.5)
$$

\vspace{3mm}

{\bf Proof.} Let us consider the following example (test function)

$$
f_0(x) = x^{-\beta} \cdot I(x \in (0,1)), \ \beta < 1/p,
$$
and denote $  g_0(x) = I^{\alpha}[f_0](x)/\Gamma(1- \alpha); $ then

$$
 \forall \beta \in (0,1/p) \Rightarrow f_0 \in L_p(0,1) \subset L_p(0, \infty):
$$

$$
|f_0|_p = (1 - \beta p)^{-1/p} < \infty.
$$

 Further, let $  x \in (0,1); $ we have consequently

$$
g_0(x) = \int_0^x \frac{y^{-\beta} \ dy}{(x - y)^{1 - \alpha}} = x^{\alpha - \beta} \cdot B(1 - \beta,\alpha)
$$
and therefore

$$
\omega(g_0,h) = C(\alpha,\beta) \cdot h^{\alpha - \beta}, \ C(\alpha,\beta) =  B(1 - \beta,\alpha) > 0,  \ h \in (0,1). \eqno(3.6)
$$
Since the value $ \beta  $ is arbitrary from the set $ (0, 1/p), $ the proposition is proved.\par

\vspace{3mm}

{\bf B. \ Case } $ X = \{x\} = (0,\infty).  $ \par

\vspace{3mm}

{\bf Proposition 3.B.} Suppose the inequality

$$
\omega \left(I^{\alpha}([f],h) \right) \le  F_{\alpha,p}(f) \cdot h^{\mu(\alpha,p)}, \ p > 1/\alpha. \eqno(3.7)
$$
there holds for any $ f \in L_p(R_+)  $   Then

$$
\mu(\alpha,p) = \alpha - 1/p. \eqno(3.8)
$$

\vspace{3mm}

{\bf Proof.} We will use the well - known {\it scaling,}  or {\it  dilation} method, see \cite{Stein1}, chapter 3,
\cite{Samko1}, chapter 3, \cite{Talenty1}. Namely, let $ \rho(\cdot)  $ be arbitrary non - zero function from the space $  L_p(R_+) $
such that

$$
\omega(I^{\alpha}[\rho],h) \le K \ h^{\gamma} \ |\rho|_p, \ h > 0. \eqno(3.9)
$$
 Let also $ \lambda  $ be arbitrary positive number; the dilation (linear) operator $  T_{\lambda} $ is defined by an equality

$$
T_{\lambda}[\rho](x) \stackrel{def}{=} \rho(\lambda x).
$$
 Evidently, $   T_{\lambda}[\rho](\cdot) \in L_p(R_+) $ and moreover

$$
|T_{\lambda} [\rho]|_p = \lambda^{-1/p} |\rho|_p
$$
and analogously

$$
I^{\alpha} T_{\lambda}[\rho] = \lambda^{-\alpha} T_{\lambda}I^{\alpha}[\rho].  \eqno(3.10)
$$

 We deduce from the source inequality (3.8) applied to the function $ T_{\lambda}[\rho]  $

$$
\lambda^{-\alpha} \omega(I^{\alpha} \rho, \lambda h) \le K \ \lambda ^{-1/p} \ h^{\gamma} \ |\rho|_p,
$$
or equally after change of variables $  \lambda h = \delta, \ \delta \in (0,\infty) $

$$
\omega(I^{\alpha}[\rho], \delta) \le K \ \delta^{\gamma} \ \lambda^{\alpha - \gamma - 1/p} \ |\rho|_p. \eqno(3.11)
$$
Thus, $ \alpha - \gamma - 1/p = 0,  $ Q.E.D.\par

\vspace{4mm}

 \section{Module of continuity of Riesz potential.}

  The necessary and sufficient condition for existence (a.e.) of Riesz potential (1.3) is the following

$$
\int_{R^d} (1 + |y|)^{\alpha - d} \ |f(y)| \ dy < \infty; \eqno(4.1)
$$
we will suppose in what follows in this and in the next sections  that this condition on the (measurable) function $  f: R^d \to R  $
is satisfied. \par
 Further, assume $  \alpha \in (0,d), \ f \in L_p(R^d) $ for some $  p > d/\alpha: $

$$
\int_{R^d} |f(y)|^p  \ dy < \infty, \ p = \const > d/\alpha. \eqno(4.2)
$$

 We introduce therefore the following norm for a (measurable) function $ f: R^d \to R $

$$
|f|_{\alpha,d,p} \stackrel{def}{=} \max \left\{\int_{R^d} (1 + |y|)^{\alpha - d} \ |f(y)| \ dy, \  |f|_p \right\}  \eqno(4.3)
$$
and the correspondent Banach  space consisting on all the measurable functions $ F: R^d \to R  $ with finite norm

$$
L_{\alpha,d,p} = \left\{ f, f: R^d \to R,  \ |f|_{\alpha,d,p} < \infty  \right\}.
$$

 Y.Mizuta et all in \cite{Mizuta1}-\cite{Mizuta4} proved that under condition $  |f|_{\alpha,d,p} < \infty $
the Riesz potential

 $$
 r_{\alpha}[f](x):= R_{\alpha}[f](x) \eqno(4.4)
 $$
is uniformly continuous and moreover

$$
\omega(r_{\alpha}[f], h) \le C \cdot  \left[ \frac{p-1}{\alpha p - d} \right]^{1 - 1/p} \cdot h^{\alpha - d/p} \cdot |f|_{\alpha,d,p}, \ h > 0.
\eqno(4.5)
$$

 It turns out that the exponent $ \alpha - d/p $ in the considered case is non - improvable still for the Riesz potential. \par

We denote by $  K_R(\alpha,p) $  the optimal, i.e. minimal value $  C  $ in the last inequality:

$$
K_R(\alpha,p) \stackrel{def}{=} \sup_{h > 0} \sup_{0 \ne f \in L_{\alpha,d,p}}
\left[ \frac{\omega(r_{\alpha}[f], h) \cdot (\alpha p - d)^{1 - 1/p} }{ (p-1)^{1 - 1/p} \cdot h^{\alpha - d/p} \cdot |f|_{\alpha,d,p }} \right],
\eqno(4.6)
$$
then as before for some finite positive constants  $    C_1(d), \  C_2(d) $ depending only on the dimension $  d  $

$$
\frac{C_1(d)}{\alpha} \le K_R(\alpha,p) \le \frac{C_2(d)}{\alpha (d - \alpha)}, \ 0 < \alpha < d.
$$

\vspace{3mm}

{\bf Proposition 4.C.} Suppose the inequality

$$
\omega \left(R_{\alpha}([f],h) \right) \le  F_{\alpha,p}(f) \cdot h^{\nu(\alpha,p)}, \ p > d/\alpha, \ h \ge 0
$$
there holds for arbitrary $ f \in L_{\alpha,d,p}(R^d),  $ i.e.
and such that  $ \int_{R^d} (1 + |y|)^{\alpha - d} \ |f(y)| \ dy < \infty $ and $ f \in L_p(R^d) $  Then

$$
\nu(\alpha,p) = \alpha - d/p. \eqno(4.7)
$$

\vspace{3mm}

{\bf Proof} is at the same as in the proposition 3.B  by means of scaling method an may be omitted.\par

\vspace{3mm}

 Further, let us introduce the following Young-Orlicz function

$$
\Phi_{p,\gamma}(u) = |u|^p \ (\ln |u|)^{\gamma}, \ |u|>e, \eqno(4.8a)
$$

$$
\Phi_{p,\gamma}(u) = e^{p-2} \ u^2, \ |u| < e; \ p = \const > 1, \ \gamma = \const > 0; \eqno(4.8b)
$$
and denote the correspondent Orlicz space by $ L( \Phi_{p,\gamma}) $  with a norm $ ||f||L( \Phi_{p,\gamma}) . $
 Since this function $ u \to \Phi_{p,\gamma}(u) $ satisfies the $  \Delta_2 $
condition, the belonging of arbitrary (measurable) function $  f: R^d \to R  $ to this  Orlicz space: $  f \in L( \Phi_{p,\gamma})  $
is completely equivalent to the convergence of the following integral

$$
\int_{R^d} \Phi_{p,\gamma} (f(x)) \ dx < \infty. \eqno(4.9)
$$

 Introduce also the following norm

$$
 ||f||L_{\alpha,d}( \Phi_{p,\gamma}) \stackrel{def}{=} ||f||L( \Phi_{p,\gamma}) + \int_{R^d} (1 + |y|)^{\alpha - d} \ |f(y)| \ dy. \eqno(4.10)
$$

 It follows immediately from the articles \cite{Mizuta1}  -  \cite{Mizuta4} that if $ ||f||L_{\alpha,d}( \Phi_{p,\gamma}) < \infty, $ then

$$
 \omega \left(R_{\alpha}([f],h) \right) \le C(d) \cdot  \left[ \frac{p-1}{\alpha p - d} \right]^{1 - 1/p} \times
$$

$$
  \ h^{\alpha - d/p} \ |\ln h|^{\gamma/p} \cdot ||f||L_{\alpha,d}( \Phi_{p,\gamma}),
 \ 0 < h < 1/e. \eqno(4.11)
$$

\vspace{3mm}

 \section{Module of continuity of fractional integrals for the functions from Grand Lebesgue Spaces.}

\vspace{3mm}

  Recently  appear the so-called Grand Lebesgue Spaces $ GLS = G(\psi) =G\psi =
 G(\psi; A,B), \ A,B = \const, A \ge 1, A < B \le \infty, $ spaces consisting
 on all the measurable functions $ f: R^d \to R $ with finite norms

$$
   ||f||G(\psi) \stackrel{def}{=} \sup_{p \in (A,B)} \left[ |f|_p /\psi(p) \right].
\eqno(5.1)
$$

  Here $ \psi(\cdot) $ is some continuous positive on the {\it open} interval
$ (A,B) $ function such that

$$
     \inf_{p \in (A,B)} \psi(p) > 0, \ \psi(p) = \infty, \ p \notin (A,B).
$$
 We will denote
$$
 \supp (\psi) \stackrel{def}{=} (A,B) = \{p: \psi(p) < \infty, \} \eqno(5.2)
$$

The set of all $ \psi $  functions with support $ \supp (\psi)= (A,B) $ will be
denoted by $ \Psi(A,B). $ \par
  This spaces are rearrangement invariant, see \cite{Bennet1}, and
    are used, for example, in the theory of probability  \cite{Kozachenko1},
  \cite{Ostrovsky1}, \cite{Ostrovsky2}; theory of Partial Differential Equations
  \cite{Iwaniec2};  functional analysis \cite{Iwaniec2},
  \cite{Ostrovsky2}; theory of Fourier series, theory of martingales, mathematical statistics,
  theory of approximation etc.\par

 Notice that in the case when $ \psi(\cdot) \in \Psi(A,\infty)  $ and a function
 $ p \to p \cdot \log \psi(p) $ is convex,  then the space
$ G\psi $ coincides with some {\it exponential} Orlicz space. \par
 Conversely, if $ B < \infty, $ then the space $ G\psi(A,B) $ does  not coincides with
 the classical rearrangement invariant spaces: Orlicz, Lorentz, Marcinkiewicz  etc.\par

  The fundamental function of these spaces $ \phi(G(\psi), \delta) = ||I_A ||G(\psi), \mes(A) = \delta, \ \delta > 0, $
where  $ I_A  $ denotes as ordinary the indicator function of the measurable set $ A, $ by the formulae

$$
\phi(G(\psi), \delta) = \sup_{ p \in \supp (\psi)} \left[ \frac{\delta^{1/p}}{\psi(p)} \right].
\eqno(5.2)
$$
The fundamental function of arbitrary rearrangement invariant spaces plays very important role in functional analysis,
theory of Fourier series and transform \cite{Bennet1} as well as in our further narration. \par

 Many examples of fundamental functions for some $ G\psi $ spaces are calculated in  \cite{Ostrovsky1}, \cite{Ostrovsky2}.\par

\vspace{3mm}

{\bf Remark 5.1.} If we introduce the {\it discontinuous} function

$$
\psi_{(r)}(p) = 1, \ p = r; \psi_{(r)}(p) = \infty, \ p \ne r, \ p,r \in (A,B)
$$
and define formally  $ C/\infty = 0, \ C = \const \in R^1, $ then  the norm
in the space $ G(\psi_r) $ coincides with the $ L_r $ norm:

$$
||f||G(\psi_{(r)}) = |f|_r.
$$
Thus, the Grand Lebesgue Spaces are direct generalization of the
classical exponential Orlicz's spaces and Lebesgue spaces $ L_r. $ \par

\vspace{3mm}

 Suppose that the function $  f: R^d \to R $ and parameters $ \alpha, d,p $ satisfy  all the conditions of the
previous section. Suppose also that the function $  p \to |f|_{\alpha,d,p} $ allows the following estimation

$$
  |f|_{\alpha,d,p} \le \psi_{\alpha,d}(p), \eqno(5.3)
$$
 where

$$
\psi_{\alpha,d}(\cdot) \in \Psi(A,B), \ \exists A = \const > d/\alpha, \ \exists B > A. \eqno(5.3a)
$$

 One can choose, for instance,

$$
\psi_{\alpha,d}(p) :=  |f|_{\alpha,d,p}, \eqno(5.3b)
$$
if the function  $ p \to \psi_{\alpha,d}(p) $ satisfies of the condition (5.3a). \par

 Define a new $  \psi \ - $ function $  \nu_{\alpha,d}(p) $ as follows:

$$
\nu_{\alpha,d}(p) := \psi_{\alpha,d}(p)\cdot K_R(\alpha,p) \cdot
\left[ \frac{p-1}{\alpha p - d} \right]^{1 - 1/p}, \ p \in (A,B). \eqno(5.4)
$$

\vspace{3mm}

{\bf Theorem 5.1.} We propose under formulated above conditions

$$
\omega \left( r_{\alpha}[f], \delta \right) \le   \frac{\delta^{\alpha}}{\phi(G\nu_{\alpha,d}(p), \delta^d) }, \ \delta > 0. \eqno(5.5)
$$

\vspace{3mm}

{\bf Proof.} We use the Mizuta's inequality and direct definition of the constant $ K_R(\alpha,p)  $

$$
\omega(r_{\alpha}[f], h) \le  K_R(\alpha,p) \cdot  \left[ \frac{p-1}{\alpha p - d} \right]^{1 - 1/p} \cdot h^{\alpha - d/p} \cdot |f|_{\alpha,d,p} \le
$$

$$
K_R(\alpha,p) \ \psi_{\alpha,d,p}(p) \ h^{\alpha - d/p} \le \nu_{\alpha,d,p}(p) \ h^{\alpha - d/p}, \ p \in (A,B).
$$
Therefore

$$
\frac{\omega(r_{\alpha}[f], h)}{ h^{\alpha}} \le \frac{h^{-d/p}}{1/\nu_{\alpha,d,p}(p)} =
\left[ \frac{h^{d/p}}{ \nu_{\alpha,d,p}(p)} \right]^{-1}. \eqno(5.6)
$$

 Since the left - hand  side of the last inequality does not dependent on the variable $ p, $ we can take the infinum from
both the sides one:

$$
\frac{\omega(r_{\alpha}[f], h)}{ h^{\alpha}} \le
\inf_{p \in (A,B)} \left\{ \frac{h^{-d/p}}{1/\nu_{\alpha,d,p}(p)} \right\} =
 \left[ \sup_{p \in (A,B)}  \left\{ \frac{h^{d/p}}{ \nu_{\alpha,d,p}(p)} \right\} \right]^{-1} =
$$

$$
\frac{1}{\phi(G\nu_{\alpha,d,p}, h^d)}, \eqno(5.7)
$$
which is equivalent to the proposition of theorem. \par

\vspace{3mm}

{\it Let us consider a slight generalization of theorem 5.1.} \par

\vspace{3mm}

 Suppose $ ||f||L_{\alpha,d}( \Phi_{p,\gamma}) < \infty. $
Suppose also that the function $  p \to ||f||L_{\alpha,d}( \Phi_{p,\gamma}) < \infty $ allows the following estimation

$$
||f||L_{\alpha,d}( \Phi_{p,\gamma}) \le \theta_{\alpha,d}(p), \eqno(5.8)
$$
 where

$$
\theta_{\alpha,\gamma, d}(\cdot) \in \Psi(A_1,B_1), \ \exists A_1 = \const > d/\alpha, \ \exists B_1 > A_1. \eqno(5.8a)
$$

 One can choose, for instance,

$$
\theta_{\alpha,\gamma, d}(p):= ||f||L_{\alpha,d}(\Phi_{p,\gamma}),
$$
if the function  $ p \to ||f||L_{\alpha,d}(\Phi_{p,\gamma}) $ satisfies of the condition (5.8a). \par

 Define a new $  \psi \ - $ function $  \zeta_{\alpha,\gamma, d}(p) $ as follows:

$$
\zeta_{\alpha,\gamma, d}(p) := \theta_{\alpha, \gamma, d}(p)\cdot K_R(\alpha,p) \cdot
\left[ \frac{p-1}{\alpha p - d} \right]^{1 - 1/p}, \ p \in (A_1,B_1). \eqno(5.10)
$$

\vspace{3mm}

{\bf Theorem 5.2.} We propose under formulated above conditions

$$
\omega \left( r_{\alpha}[f], h \right) \le
\frac{h^{\alpha} }{\phi(G\zeta_{\alpha,\gamma, d}, h^d  \cdot |\ln h|^{-\gamma}) }, \ h \in (0, 1/e). \eqno(5.11)
$$

\vspace{3mm}

Consider ultimately the case when $ \exists (A_2, B_2) = \const, \ A_2 > d/\alpha, B_2 \in (A_2, \infty],  $ such that

$$
 \forall p \in (A_2, B_2) \ \Rightarrow  \int_{R^d} |f(y)|^p \ [ \ln_+ |f(y)|]^{\gamma_0 + \gamma_1 p} \ dy < \infty, \eqno(5.12)
$$
where $ \ln_+ z = \max(1, \ln z), \ z \ge 0, $ for certain finite non-negative constants $ \gamma_0, \ \gamma_1. $ We introduce a
non-negative function $  \kappa_0 = \kappa_0(p) $ by an equality

$$
\kappa^p_0(p) := \int_{R^d} |f(y)|^p \ [\ln_+|f(y)|]^{\gamma p} \ dy =
\left| f \cdot \left[\ln_+|f(y)|]^{\gamma}\right] \right|_p^p,  \ p \in (A_2, B_2),
$$
and define the new finite in the interval $  (A_2, B_2) $ function

$$
\kappa(p) := \max \left\{\int_{R^d} (1 + |y|)^{\alpha - d} \ |f(y)| \ dy, \ \kappa_0(p) \right\}. \eqno(5.13)
$$

\vspace{3mm}

{\bf Theorem 5.3.} We assert under formulated above conditions, for instance, conditions (4.1) and (5.12)

$$
\omega \left( r_{\alpha}[f], h \right) \le
\frac{C(\alpha,\gamma_0,\gamma_1,d) \ h^{\alpha}\ |\ln h|^{-\gamma_1} }{\phi(G\kappa, h^d
 \cdot |\ln h|^{-\gamma_0}) }, \ h \in (0, 1/e). \eqno(5.14)
$$

 \vspace{4mm}

 \section{Concluding remarks.}

 \vspace{3mm}

 It is interest perhaps to obtain the estimations of the module continuity
for the {\it weight} Riesz potential as well as for the {\it weight} fractional integrals and derivatives.
The $ L_p $ estimates for ones  are investigates in \cite{Lieb1}, \cite{Lieb2},
\cite{Lieb3}, \cite{Ostrovsky12}, \cite{Ostrovsky13}.\par
 Also it is interest by our opinion to investigate the multidimensional fractional integrals and derivatives.\par

\end{document}